\documentclass[11pt]{amsart}
\usepackage[english]{babel}

\usepackage{amssymb,amsmath,amscd,amsthm}

\def\{{\protect\lbrace}
\def\}{\protect\rbrace}

\newcommand{\Kdim}{\operatorname{Kdim}}

\begin{document}

\begin{center}
\textbf{\large Arithmetical Rings and Krull Dimension}
\end{center}

\hfill {\sf A.A. Tuganbaev}

\hfill National Research University "MPEI", Moscow, Lomonosov Moscow State University

\hfill e-mail: tuganbaev@gmail.com

{\bf Abstract.} Let $A$ be a commutative arithmetical ring. The ring $A$ has Krull dimension if and only if every factor ring of $A$ is finite-dimensional and does not have idempotent proper essential ideals.\\
The study is supported by Russian Science Foundation (project no.~16-11-10013).


All rings are assumed to be associative and with zero identity element; all modules are unitary right modules. A ring $A$ is said to be {\it arithmetical} if the lattice of two-sided ideals of $A$ is distributive. A ring is said to be \textit{right} (resp., \textit{left}) \textit{uniserial} if any two of its right (resp., left) ideals are comparable with respect to inclusion. For example, every residue ring $\mathbb{Z}/n\mathbb{Z}$ is a finite commutative arithmetical ring. This ring is uniserial if and only if $n$ is a power of a prime integer. This paper is a continuation of the paper \cite{Tug15} devoted to arithmetical rings.

 We recall the transfinite definition of the Krull dimension $\Kdim M$ of a module $M$, see\cite{GorR73}. (We note that not all modules have Krull dimension.) By definition, we assume that zero modules have Krull dimension, which is equal to $-1$, and every non-zero Artinian module has Krull dimension which is equal to zero. Let us assume that $\alpha >0$ is an ordinal, the modules with Krull dimension $\beta$ are defined for all all ordinal numbers $\beta <\alpha$, and $M$ is a module such that $\Kdim M\ne \beta$. One says that the {\it Krull dimension} $\Kdim M$ is equal to $\alpha$ if for any infinite properly descending chain $M_1>M_2>\ldots$ of submodules in $M$, there exists a positive integer $n$ such that $\Kdim (M_n/M_{n+1})<\alpha$. All factor modules and submodules of a module with Krull dimension have Krull dimension, every extension of a module with Krull dimension by a module with Krull dimension has Krull dimension, and every Noetherian module has Krull dimension, see \cite{GorR73}. For a ring $A$, the \textit{right Krull dimension} $\Kdim (A_A)$ is the Krull dimension of the module $A_A$ if it exists. For example, the Krull dimension of the residue ring $\mathbb{Z}/n\mathbb{Z}$ is equal to $0$, the Krull dimension of the polynomial ring $k[x_1,\ldots ,x_n]$ over any field $k$ is equal to $n$. (For a commutative ring $A$, the Krull dimension of $A$ is equal to the classical Krull dimension defined with the use of prime ideals.) The class of all rings with right Krull dimension is larger than the class of all right Noetherian rings and rings with right Krull dimension satisfy many useful properties of right Noetherian rings. For example, if $A$ is a ring with right Krull dimension, then its prime radical $P$ is nilpotent and the factor ring $A/P$ has the right classical ring of fractions which is a semisimple Artinian ring.

A module $M$ is said to be \textit{finite-dimensional} if $M$ does not contain a direct sum of infinite number of non-zero submodules. Each module with Krull dimension is finite-dimensional, see \cite[Proposition 1.4]{GorR73}. If $M$ is a module and $X$ is a submodule in $M$ such that $X\cap Y\ne 0$ for any non-zero submodule $Y$ in $M$, then $X$ is called an \textit{essential} submodule in $M$. An ideal $B$ is said to be \textit{idempotent} if $B=B^2$. 

\textbf{Remark 1.} In \cite{Lem77}, it is proved that a commutative uniserial ring $A$ has Krull dimension if and only if $A$ does not have idempotent proper non-zero ideals. In this case, it is clear that every factor ring of the ring $A$ does not have idempotent proper non-zero ideals. It is clear that every factor ring of any commutative uniserial ring is a finite-dimensional uniserial ring such that every non-zero ideal in $A$ is essential. We remark that $\mathbb{Z}$ is a commutative arithmetical ring with Krull dimension $1$, the ring $\mathbb{Z}$ does not have idempotent proper non-zero ideals and the factor ring $\mathbb{Z}/6\mathbb{Z}$ of $\mathbb{Z}$ has an idempotent proper non-zero ideal.

\textbf{Remark 2.} The direct product $A=\prod _{i\in I}A_i$ of any set of fields $A_i$ is a commutative arithmetical ring such that all ideals of $A$ are idempotent ideals. The ring $A$ has Krull dimension if and only if the set $I$ is finite.

\textbf{Remark 3.} Let $A$ be a commutative ring with Krull dimension. In \cite[Th\'{e}or\`{e}me 12]{Lem72}, it is proved that $A=A_1\times\ldots \cdot \times A_n$, where each of the rings $A_i$ has Krull dimension and does not have idempotent proper non-zero ideals. Since every factor ring of the ring $A$ has Krull dimension, we have that every factor ring $\overline{A}$ of the ring $A$ does not have idempotent proper essential ideals. In addition, $\overline{A}$ is a finite-dimensional ring.

In connection to Remarks 1, 2, and 3, we prove Theorem 1 which is the main result of this paper.

\textbf{Theorem 1.} \textit{Let $A$ be a commutative arithmetical ring. The following conditions are equivalent.
\begin{enumerate}
\item[\textbf{1)}] 
$A$ has Krull dimension.
\item[\textbf{2)}] 
Every factor ring of the ring $A$ is finite-dimensional and does not have idempotent proper essential ideals.
\end{enumerate}}

The proof of Theorem 1 is decomposed into a series of assertions, some of which are of independent interest. 

A module $M$ is said to be \textit{distributive} if the lattice of all submodules of $M$ is distributive, i.e. $X\cap (Y+Z)=X\cap Y+X\cap Z$ for any two submodules $X,Y,Z$ of the module $M$. A module $M$ is said to be \textit{uniform} if any two non-zero submodules of $M$ have the non-zero intersection. A ring $A$ is called a \textit{domain} (resp., a \textit{prime} ring) if the product of any two non-zero elements (resp., ideals) of $A$ is non-zero. For a ring $A$, we denote by $J(A)$ the Jacobson radical of $A$.  A proper ideal $P$ of the ring $A$ is said to be \textit{completely prime} (resp., \textit{prime}) if the factor ring $A/P$ is a domain (resp., a prime ring). A right ideal $P$ of the ring $A$ is said to be \textit{completely prime} if $ab\notin P$ for all $a,b\in A\setminus P$. The intersection of all prime ideals of the ring $A$ is a nil-ideal and is called the \textit{prime radical} of $A$. A ring $A$ is said to be \textit{right invariant} (resp.,  \textit{left invariant}) if all right (resp., left) ideals of $A$ are ideals. A ring $Q$ is called the \textit{classical right ring of fractions} of the ring $A$ if $A$ is a unitary subring in $Q$, every non-zero-divisor of the ring $A$ is invertible in $Q$, and for each element $q\in Q$, there exist elements $a,s\in A$ such that $s$ is a non-zero-divisor in $A$ and $q=as^{-1}$.

\textbf{Lemma 1 \cite{Ste74}.} \textit{Let $A$ be a ring and let $X$ be a right $A$-module. The module $X$ is distributive if and only if for any two elements $x,y\in M$, there exist elements $a,b\in A$ such that $a+b=1$ and $xaA+ybA\subseteq xA\cap yA$.}

\textbf{Lemma 2.} \textit{Let $A$ be a right distributive ring.
\begin{enumerate}
\item[{\bf 1)}] 
If $P_1\subset P_2\subset \ldots $ is an infinite properly ascending chain of completely prime right ideals of the ring $A$, then $\cup_{i=1}^{\infty}P_i$ is an idempotent proper completely prime right ideal.
\item[{\bf 2)}] 
If the ring $A$ does not have idempotent proper completely prime ideals, then $A$ is a ring with the maximum conditions on completely prime ideals.
\item[{\bf 3)}] 
If the ring $A$ does not have idempotent proper completely prime right ideals, then $A$ is a ring with with the maximum conditions on completely prime right ideals.
\end{enumerate}}

\textbf{Proof.} \textbf{1).} We denote by $X$ the right ideal $\cup_{i=1}^{\infty}P_i$. Since all $P_i$ are proper right ideals, $X$ is a proper right ideal. Since all $P_i$ are completely prime right ideals, $X$ is a completely prime right ideal.

We assume that $X\ne X^2$. Let $x\in X\setminus X^2$. Then $x\in P_i$ for some $i$. Since $P_i$ is a completely prime right ideal and $X$ properly contains $P_i$, we have that $X^2\not\subseteq P_i$ and there exists an element $y\in X^2\setminus P_i$. By Lemma 1, there exist elements $a,b\in A$ such that $a+b=1$, $xa\in yA\subseteq X^2$ and $yb\in xA\subseteq P_i$. Since $x\notin X^2$ and $xa\in X^2$, we have $xb=x-xa\notin X^2$. Therefore, $b\notin X$. Since $y\notin P_i$ and the element $yb$ is contained in the completely prime right ideal $P_i$, we have $b\in P_i\subseteq X$. This is a contradiction.

\textbf{2) and 3).} The assertions follow from 1).~\hfill$\square$

\textbf{Lemma 3.} \textit{Let $A$ be a right invariant arithmetical ring which does not have idempotent proper prime ideals. Then $A$ is a ring with the maximum conditions on prime ideals.}

\textbf{Proof.} Since $A$ is a right invariant arithmetical ring, $A$ is a right distributive ring. In addition, every prime ideal of the right invariant ring $A$ is completely prime. Therefore, the assertion follows from Lemma 2(2).~\hfill$\square$

\textbf{Lemma 4 \cite[Lemma 1.6]{Tug02}.} \textit{Let $A$ be a right distributive ring such that the set $N$ of all right or left zero-divisors of the ring $A$ is a completely prime ideal of the ring $A$. Then $A$ has the right classical right ring of fractions $Q$, which is a right uniserial ring, and $NQ=J(Q)$.}

\textbf{Lemma 5.} \textit{Let $A$ be a right distributive right uniform ring such that all right zero-divisors of $A$ are left zero-divisors and let $N$ be the set of all zero-divisors of the ring $A$ . Then $N$ is a completely prime ideal of $A$, the ring $A$ has the right uniserial classical right ring of fractions $Q$ and $NQ=J(Q)$.\\
In addition, if $ A $ has a completely prime nil-ideal $P$, then $ P=xP $ for any non-zero-divisor $x$ of the ring $A$, $ P=QP $ is a left ideal of the ring $Q$, $ PQ=QPQ $ is an ideal of the ring $Q$, and $ (PQ)^n=P^nQ $ for every positive integer $ n $.}

\textbf{Proof.} Let $N_1$ be the set of all left zero-divisors in $A$. By \cite[Lemma 1.4]{Tug02}, $N_1$ is a completely prime ideal of the ring $A$. Since all right zero-divisors in the ring $A$ are left zero-divisors, $N=N_1$. Therefore, $N$ is a completely prime ideal. By Lemma 4, $A$ has the right classical right ring of fractions $Q$, which is a right uniserial ring, and $NQ=J(Q)$.

We assume that $ A $ has a completely prime nil-ideal $P$. Let $y\in P$. By Lemma 1, there exist elements $a,b\in A$ such that $a+b=1$ and $xa,yb\in xA\cap yA$. Since $ x\notin P $ and $xa$ is an element of the completely prime ideal $P$, we have that $a$ is an element of the nil-ideal $P$. Then $b=1-a$ is an invertible element, and $yb\in xA$. Therefore, $y=ybb^{-1}\in xA$ and $P=xP$. Since $x$ is an invertible element of the ring $Q$ and $P=xP$, we have $P=x^{-1}P$. Therefore, $P=QP$. Then $ PQ=QPQ $ is an ideal of the ring $ A $ and $ (PQ)^n=P^nQ $ for every positive integer $ n $.~\hfill$\square$

\textbf{Lemma 6.} \textit{Let $A$ be a commutative arithmetical uniform ring, $P$ be the prime radical of the ring $A$, and $A/P$ be a finite-dimensional ring. 
\begin{enumerate}
\item[\textbf{1)}] 
$P$ is a completely prime nil-ideal and either $P=P^2$, or $P$ is a nilpotent ideal.
\item[\textbf{2)}] 
If the ring $A$ does not have idempotent proper essential ideals, then $P$ is a nilpotent ideal.
\end{enumerate}}

\textbf{Proof.} \textbf{1).} Since $A/P$ is a commutative finite-dimensional semiprime arithmetical ring, $A/P$ is a finite direct product of domains \cite[Proposition 2]{Cam75}. Since $P$ is a nil-ideal, all idempotents of the ring $A/P$ can be lifted to idempotents of the commutative indecomposable ring $A$. Therefore, the ring $A/P$ does not have non-trivial idempotents. Then $A/P$ is a domain and $P$ is a completely prime nil-ideal. By Lemma 5, the ring $A$ has the uniserial classical ring of fractions $Q$, $ P=xP $ for any non-zero-divisor $x$ of the ring $A$, $ P=QP $ is an ideal of the ring $Q$, and $P^n=P^nQ$ for every positive integer $ n $. 

We assume that $P\ne P^2$. Let $p\in P\setminus P^2$. Then
 $pQ$ is a nilpotent ideal of the commutative uniserial ring $Q$ which is not contained in the ideal $P^2=P^2Q$. Then ideal $P^2$ is contained in the nilpotent ideal $pQ$. Therefore, the ideal $P^2$ is a nilpotent ideal. Therefore, $P$ is a nilpotent ideal.

\textbf{2).} Without loss of generality, we can assume that $P$ is a non-zero proper ideal of the ring $A$. Since $A$ is a uniform ring, $P$ is an essential ideal. By assumption, the ring $A$ does not have proper essential idempotent ideals. Therefore, $P\ne P^2$. By 1), $P$ is a nilpotent ideal.~\hfill$\square$

\textbf{Lemma 7 \cite[Proposition 2]{Lem77}.} \textit{Let $M$ be a module such that all factor modules of the module $M$ are finite-dimensional and each non-zero factor module $Q$ of the module $M$ contains a non-zero submodule with Krull dimension. Then $M$ has Krull dimension.}

\textbf{Lemma 8 \cite[1.4, 7.1, 7.3, 7.4]{GorR73}.} \textit{Let $A$ be a ring with right Krull dimension. Then every factor ring of the ring $A$ is right finite-dimensional, $A$ is a ring with the maximum condition on prime ideals, and for every proper ideal $B$ of the ring $A$, there exist prime ideals $P_1,\ldots , P_n$ of the ring $A$ such that $P_1\cdot\ldots \cdot P_n\subseteq B$ and each of the ideals $P_i$ contains the ideal $B$.}

The proofs of the following Lemma 9 and Lemma 10 are based on the ideas of the paper \cite{Lem77}.

\textbf{Lemma 9.} \textit{Let $A$ be a ring such that all cyclic right $A$-modules are finite-dimensional.}
\begin{enumerate}
\item[\textbf{1)}] 
\textit{If there exists a finite set $\{P_1,\ldots ,P_n\}$ ideals of the ring $A$ such that each of the cyclic right $A$-modules $A/P_i$ has Krull dimension, then the ring $A/(P_1\cdot\ldots \cdot P_n)$ has right Krull dimension.}
\item[\textbf{2)}] 
\textit{We assume that there exists a set $\mathcal{Q}$ of proper ideals of the ring $A$ such that all ascending ideals from $\mathcal{Q}$ stabilize at finite step and for any non-zero right ideal $B$ of the ring $A$, there exist ideals $P_1,\ldots,P_n\in\mathcal{Q}$ such that
$$
P_1\cdot \ldots \cdot P_n\subseteq B \subseteq P_1\cap\ldots\cap P_n.
$$
Then $A$ has right Krull dimension.}
\end{enumerate}

\textbf{Proof.} \textbf{1).} The proof of this assertion is reported to the author by V.T. Markov. We use the induction on $n$. For $n=1$, it is nothing to prove. Let $n>1$ and let $Q=P_1\ldots P_{n-1}$. By the induction assumption, the factor ring $A/Q$ has right Krull dimension. Let $X$ be an arbitrary non-zero factor module of the right $A$-of the module $Q/QP_n$. Then $X$ is a submodule of the right $A/P_n$-module $Q/QP_n$. Since $X$ is a right module over the ring $A/P_n$ with right Krull dimension, every non-zero submodule of the $A/P_n$-module $X$ contains a non-zero submodule with Krull dimension. Then every non-zero submodule of the $A$-module $X$ contains a non-zero submodule with Krull dimension. In addition, $X_A$ is a submodule of the cyclic $A$-module $A/QP_n$ and all factor modules of the cyclic $A$-module $A/QP_n$ are finite-dimensional by the assumption. Therefore, all factor modules of the module $X_A$ are finite-dimensional. By Lemma 7, the module $X_A$ has Krull dimension. Since $X$ is an arbitrary non-zero factor module $(Q/QP_n)_A$, the module $(Q/QP_n)_A$ has Krull dimension. In addition, $(A/Q)_A$ has Krull dimension by the induction hypothesis. Since $(A/QP_n)_A$ is an extension of the module $(Q/QP_n)_A$ with Krull dimension by the module $(A/Q)_n$, we have that $(A/QP_n)_A$ has Krull dimension Therefore, the ring $A/QP_n$ has right Krull dimension.

\textbf{2).} We assume that $A$ does not have right Krull dimension. We construct a properly ascending sequence of ideals $Q_0\subset Q_1\subset Q_2\ldots$ from $\mathcal{Q}\cup\{0\}$ such that the ring $A/Q_i$ does not have Krull dimension for all $i\geq 0$. We set
$Q_0=0$. We assume that the ideals $Q_0\subset\ldots\subset Q_k$ are
constructed, where $k\geq 0$. Since the right $A$-module $A/Q_k$ does not have
Krull dimension, there exists a right ideal $B$ of the ring $A$ such that
$B\supseteq Q_k$, $B/Q_k\neq 0$ and $(A/B)_A$ does not have
Krull dimension. By the assumption, there exist ideals $P_1,\ldots,P_n\in\mathcal{Q}$ such that
$$
P_1\cdot \ldots \cdot P_n\subseteq B \subseteq P_1\cap\ldots\cap P_n.
$$
By Lemma 7, there exists an integer $i\in\{1,\ldots,n\}$ such that the cyclic right $A$-module $A/P_i$ does not have Krull dimension. We set $Q_{k+1}=P_i$. We have an infinite ascending sequence chain $Q_1\subset Q_2\subset \ldots$ of ideals from $\mathcal{Q}$; this contradicts to the assumption.~\hfill$\square$

\textbf{Lemma 10.} \textit{For a right invariant ring $A$, the following conditions are equivalent.
\begin{enumerate}
\item[\textbf{1)}] 
The ring $A$ has right Krull dimension.
\item[\textbf{2)}] 
$A$ is a ring with the maximum condition on prime ideals, all factor of the ring of the ring $A$ are right finite-dimensional, and for any non-zero ideal $B$ of the ring $A$, there exist ideals $P_1,\ldots,P_n\in\mathcal{Q}$ such that
$$
P_1\cdot \ldots \cdot P_n\subseteq B \subseteq P_1\cap\ldots\cap P_n.
$$
\end{enumerate}}

\textbf{Proof.} The implication 1)\,$\Rightarrow$\,2) follows from Lemma 8.

2)\,$\Rightarrow$\,1). The ring $A$ is right invariant and all factor rings of the ring $A$ are right finite-dimensional. Therefore, all cyclic right $A$-modules are finite-dimensional. We denote by $\mathcal{Q}$ the set of all prime ideals of the ring $A$. Now the assertion follows from Lemma 9(2).~\hfill$\square$

\textbf{Lemma 11.} \textit{Let $A$ be a commutative arithmetical ring such that every factor ring of the ring $A$ is finite-dimensional and does not have idempotent proper essential ideals. Then $A$ has Krull dimension.}

\textbf{Proof.} By assumption, every factor ring of the ring $A$ is finite-dimensional. By Lemma 3, $A$ is a ring with the maximum condition on prime ideals. Let $B$ be an arbitrary proper ideal of the ring $A$. 
By Lemma 10, it is sufficient to prove that there exist prime ideals $Q_1,\ldots , Q_n$ of the ring $A$ such that $Q_1\cdot\ldots \cdot Q_n\subseteq B$ and each of the ideals $Q_i$ contains the ideal $B$. 

By assumption, $A/B$ is a finite-dimensional ring. Therefore, there exist ideals $B_1,\ldots , B_k$ of the ring $A$ such that $B_1\cap \ldots \cap B_k=B$ and every factor ring $A/B_i$ is a uniform ring. Let $P_1,\ldots , P_k$ be ideals of the ring $A$ such that $B_i\subseteq P_i$ and $P_i/B_i$ is the prime radical of the ring $A/B_i$, $i=1,\ldots ,k$. By assumption, every factor ring $A/P_i$ is a uniform ring. By assumption, every factor ring $A/B_i$ is a commutative arithmetical ring without idempotent proper essential ideals. By Lemma 6, $P_i/B_i$ is a prime nilpotent ideal of the ring $A/B_i$, $i=1,\ldots ,k$. Therefore, there exist positive integers $n_1,\ldots ,n_k$ such that $P_i^{n_i}\subseteq B_i$, $i=1,\ldots ,k$. Therefore, the ideal $X=P_1^{n_1}\cdot\ldots\cdot P_k^{n_k}$ is contained in the ideal $B_1\cap\ldots\cap B_k=B$ and every prime ideal $P_i$ contains the ideal $B$. This implies the required assertion.~\hfill$\square$

\textbf{Remark 4.} We give the completion of the proof of Theorem 1. The implication 1)\,$\Rightarrow$\,2) follows from Remark 3. The implication 2)\,$\Rightarrow$\,1) is proved in Lemma 11.

The author thanks I.B.Kozhuhov and V.T.Markov for their helpful discussions.


\begin{thebibliography}{9}
 
\bibitem{Cam75} Camillo V. Distributive modules // J.~Algebra. -- 1975. -- Vol.~36. -- P.~16-25.

\bibitem{GorR73} Gordon R., Robson J.C. Krull dimension // Mem. Amer. Math. Soc. -- 1973. -- no.~133. -- P.~1-78.

\bibitem{Lem72} Lemonnier B. Sur les anneaux qui ont une d\'{e}viation // C.~R. Acad. Sc. Paris. Ser.~A. -- 1972. -- T.~275. -- P.~A357-A359.

\bibitem{Lem77} Lemonnier B. Dimension de Krull et cod\'{e}viation des anneaux semi-h\'{e}r\'{e}ditaires // C.~R. Acad. Sc. Paris. Ser.~A. -- 1977. -- T.~284. -- P.~A663-A666.

\bibitem{Ste74} Stephenson W. Modules whose lattice of submodules is
distributive // Proc. London Math. Soc. -- 1974. -- Vol.~28, no.~2. -- P.~291--310.

\bibitem{Tug02} Tuganbaev A.A. Structure of distributive rings // Sb. Math. -- 2002. -- Vol.~193, no.~5. -- P.~745-760.

\bibitem{Tug15} Tuganbaev A.A. Rings whose finitely generated right ideals are quasi-projective // Discrete Math. Appl. -- 2015. -- Vol.~25, no.~4. -- P.~245–251.

\end{thebibliography}
\end{document}